# THE POWER OF MULTIFOLDS:
## FOLDING THE ALGEBRAIC CLOSURE OF THE RATIONAL NUMBERS

BY TIMOTHY Y. CHOW AND C. KENNETH FAN

ABSTRACT. We define the *n*-parameter multifold and show how to use one-parameter multifolds to get the algebraic closure of the rational numbers.

## 1. INTRODUCTION

Robert Lang described how to quintisect an arbitrary angle using techniques of origami.

How can this be? Is it not well-known that angle quintisection is impossible using the Huzita-Hatori axioms[1]? Lang quintisected using a secret weapon: the multifold. In *Angle Quintisection* [RL1], Lang explains that quintisection is impossible if only one fold is allowed at a time. However, if one allows simultaneous folds in a single origami maneuver, quintisection can be achieved. The technique is not easy to perform, but with some fiddling, it can be pulled off.

The question naturally arises: just how large is the new set of constructible numbers if one allows multifolds?

In this paper, we show that the answer is the entire algebraic closure[2] of the rational numbers. In a sense, this means that the set of origami constructible numbers is as large as possible, although we suggest a further expansion in the last section.

The contents of this paper are an elaboration of an argument the second author presented to Robert Lang in person at the June 2005 Origami USA convention in New York City.

## 2. MULTIFOLDS

A multifold is a single origami maneuver that simultaneously involves the creation of more than one crease. As Lang pointed out, some multifolds can be performed in several steps, each step involving the creation of a single crease. Such multifolds would not expand the universe of constructible numbers and Lang called them *separable*. However, some multifolds cannot be reduced to single crease folds (which would fall under the purview of the Huzita-Hatori axioms). Lang used such a *nonseparable* multifold to quintisect an arbitrary angle. We know that his multifold is nonseparable because he was able to achieve a result provably impossible under the Huzita-Hatori axioms.

In this section, we describe a categorization of multifolds that the second author first defined in messages to the Origami List. A multifold is formed by forming a sequence of folds and then "rolling" these folds along until a certain alignment condition is satisfied.

---

[1] We assume the reader is familiar with the Huzita-Hatori axioms. For a good introduction see [RL2].
[2] We assume the reader has a basic familiarity with field theory. For information on the algebraic theory involved, see any introduction to abstract algebra, such as [SL].

When the alignment is attained, all creases are made sharp. The various configurations of the creases involved in the multifold can be parameterized by some number of variables. As the creases are rolled along, a path is traced in the parameter space. If $n$ parameters are required to parameterize the creases in the multifold, we call the multifold an n-*parameter* multifold.

The casual reader may safely skip the more precise definition of an $n$-parameter multifold that we now give. Let $F_k$ be a finite sequence of folds. Let $n_k$ be a sequence of nonnegative integers and let $n$ be the sum of the $n_k$. Assume the following:

1. The fold $F_k$ can be specified in terms of a standard alignment (in the sense of Huzita-Hatori) using any preexisting references and/or references created by prior folds in the sequence and/or any $n_k$ real numbers that parameterize variable references. (Normally, these parameters correspond to measurable features of the model, such as an angle or distance.)
2. There exists a domain $D$ in Euclidean $n$-dimensional space such that $D$ is equal to the closure of its interior and the map $F$ from $D$ to physical configurations of the model (which includes creases) obtained by performing the sequence of folds $F_k$ is continuous.
3. The model can be physically manipulated in a continuous manner within the image of $F$.

**Definition.** An n-*parameter multifold* is the image, under $F$, of a point in $D$.

In the context of origami constructions, only certain $n$-parameter multifolds $F(d)$, where $d$ is in $D$, are useful. Namely, $F(d)$ must satisfy an alignment condition involving existing references (including references created by the folds in the folding sequence) that must isolate the value of at least one parameter. When an alignment condition depends on only $m$ of the $n$ parameters to isolate the value of the desired parameters, we will also call such a multifold an m-*parameter multifold*. For example, let $F_1$ be defined by a fold parallel to a fixed edge of a square sheet of origami paper and a parameter $d$ that gives the distance of the crease from that edge. Every $d$ that specifies a crease within the origami square yields a one-parameter multifold. However, it is useless as far as origami constructions go to allow any such multifold because most do not allow for $d$ to be computed. In this example, the book fold, which is specified by an appropriate alignment condition, does give a constructive multifold.

Also note that the third condition for the folding sequence makes any unfolding in the sequence risky. If the folding sequence involves unfolding a fold and refolding along a crease that intersects the unfolded crease, it is no longer possible to manipulate the model continuously: to alter the first crease, the model would have to be refolded along a discretely new location as opposed to having the new crease achieved by "rolling" the original crease to a new position.

There is a distinct practical advantage to using one-parameter multifolds. In a one-parameter multifold, one is guaranteed to discover the alignment condition by rolling the creases so long as one is in the same connected component as the alignment condition. However, in multi-parameter multifolds, one could fiddle endlessly with the paper in a desperate search for the alignment condition. The difference is similar to trying to locate a buoy in the Pacific Ocean using a telescope from an airplane by flying about with and without the benefit of knowing on which line of longitude the buoy bobs.

The good news is that one-parameter multifolds suffice because the algebraic closure of the rational numbers can be constructed using one-parameter multifolds alone, that is:

**Theorem.** Suppose that the complex number $z$ is a root of a polynomial with rational coefficients. Then the real and imaginary parts of $z$ can be constructed using one-parameter multifolds.

We remark that Lang's angle quintisection uses a one-parameter multifold.

## 3. GETTING REAL

Using well known facts about algebraic numbers, it suffices to show that we can construct only the *real* roots of polynomials with rational coefficients. Briefly, the reason is that the algebraic closure of the rational numbers is closed under complex conjugation, addition, and multiplication. Because $\operatorname{Re} a = (a + \bar{a})/2$ and $\operatorname{Im} a = -i(a - \bar{a})/2$, both the real and imaginary parts of any algebraic number are algebraic.

In fact, given an algebraic number $z$ and a polynomial with rational coefficients $p$ with $p(z) = 0$, it is possible to construct polynomials with rational coefficients with the real and imaginary parts of $z$ among their roots. This enables us to not only construct the algebraic closure of the rational numbers in the abstract, but also enables us to *solve* for the roots of any given polynomial with rational coefficients.

## 4. THE SET UP

So, in accordance with Section 3, we are given a polynomial $p(x) = a_n x^n + a_{n-1} x^{n-1} + \ldots + a_1 x + a_0$ with rational coefficients and we need to construct its real roots.

At first, it might seem that one could simply pick $x$ arbitrarily and use standard constructions for addition and multiplication to construct $p(x)$. Then, by varying $x$, find a root by noting when the constructed $p(x)$ is 0. The problem with such an approach is that it may not be possible to continuously vary $x$ because, for instance, the construction may involve unfolding and refolding or other obstructions might arise[3]. For this reason, some care must be taken to ensure that this does not happen.

To describe the method of construction, it helps to forget about origami for a moment. Instead, let us imagine that we are working with an unlimited supply of rectangular sheets of paper that can come with any dimensions that we desire. We will show how to find a real root of *p* using these sheets of paper and later show how to turn this into an origami multifold.

---

[3] In fact, Robert Lang mentioned to the second author that just such an obstruction was preventing an approach by him and Erik Demaine (using linkages) from working (as of the June, 2005 OUSA).

When we speak of using "large" or "long" sheets, by the way, what we mean is that some finite sized sheets will work and if you try the construction and discover that you ran out of paper, it means you just have to try again with larger or longer sheets!

So take a large square sheet and make a book fold through it. Unfold and orient it so that the resulting crease is vertical. We will refer to this as the *zero* crease. Please refer to Figure 1.

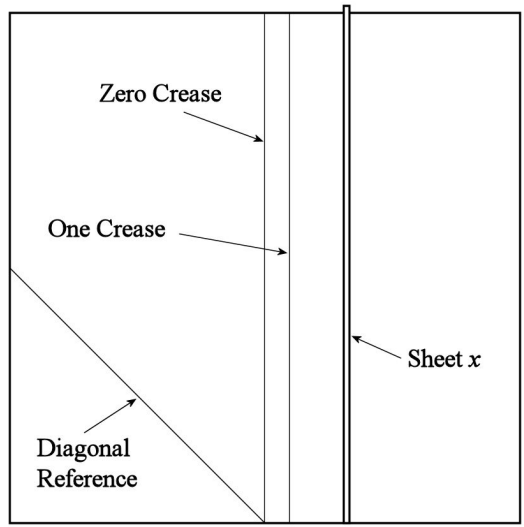

FIGURE 1. The preliminary setup.

Now, fold and unfold another sharp crease parallel to and one unit to the right of the zero crease. Call this new crease the *one* crease.

We shall assume that the polynomial is given to us by having rectangular strips with lengths equal to the absolute value of the coefficients. If you really want to, you can encode the sign of the coefficient by using one color for positive and another color for negative.

Next, fold the lower left corner of the square to the center to create a 45° crease in the lower left quadrant. Unfold. Call this crease the *diagonal reference*.

Finally, take a very long thin rectangular sheet and place it on top of the square sheet parallel to the zero and one creases. Call this *sheet x*.

Our multifold will be parameterized by the horizontal distance between the zero crease and the left edge of sheet *x*. We shall also use *x* to refer to this distance.

In the construction, we will often refer to horizontally, vertically, or diagonally aligned rectangles. Use the edges of the square or the diagonal reference to fix these orientations.

## 5. OVERVIEW OF THE CONSTRUCTION

The idea behind our construction is to define a sequence of folds that depend on the single parameter $x$. Performing the folds ultimately results in a pair of edges separated by the distance $p(x)$. This sequence gives rise to a one-parameter multifold defined by the alignment condition that this pair of edges overlap, that is, by the condition that $p(x) = 0$.

In order to do this, write $p(x)$ in the form

$$p(x) = x(x( \ldots (x(xa_n + a_{n-1}) + a_{n-2}) + \ldots ) + a_1) + a_0.$$

We devise a special sequence of folds that takes as input a pair of edges separated by the distance $d$ and outputs a pair of edges separated by the distance $xd + a$, where $a$ is a given constant. Using this special folding sequence, we can achieve a pair of edges separated by the distance $p(x)$ by working outward from the innermost embedded expression $xa_n + a_{n-1}$.

In the next section, we shall describe this special sequence of folds.

## 6. THE HEART OF THE CONSTRUCTION

Our construction will be achieved by repeated application of the following folding sequence. Please refer to Figure 2.

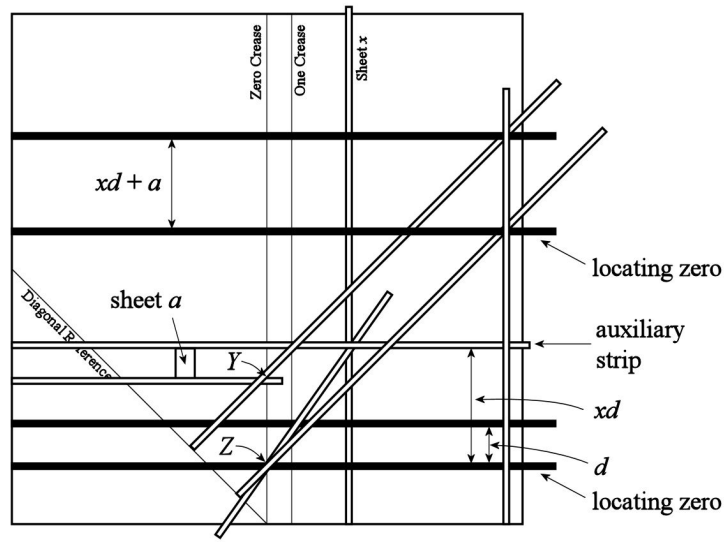

FIGURE 2. The input and output pairs of horizontal strips are solid black.

We assume that we are given a quantity $d$ encoded as the distance between the inner edges of two long, thin horizontal strips situated on the square sheet (the lower pair of black horizontal strips in Figure 2). If the quantity is positive, think of the lower strip as "locating zero", otherwise, think of the upper strip as "locating zero". When we refer to the "inner edge" of one of these long, thin horizontal strips, we shall mean "inner" in the context of the two strips forming a pair. Let $Z$ denote the point of intersection of the inner edge of the strip locating zero and the zero crease.

We aim to produce a new pair of long, thin horizontal strips that encode, in the same manner, the quantity $xd + a$ where $a$ is a real number given to us as the length of a rectangle which we will call *sheet a* (just as the coefficients of the polynomial are given to us).

To do this, note that the inner edges of the horizontal strips together with the zero and one creases form the edges of a rectangle. Place the edge of a long, thin rectangular strip along the diagonal of this rectangle so that it passes through $Z$. Follow the edge of this rectangular strip to where it intersects the left edge of sheet $x$ and mark this intersection point by placing the top or bottom (whichever is easier to see) edge of an auxiliary long, thin horizontal strip through the intersection. The vertical distance between the relevant edge of this auxiliary strip and the inner edge of the given horizontal strip locating zero is $xd$. The relevant edge of the auxiliary strip will be above or below the given horizontal strip locating zero depending on whether $xd > 0$ or $xd < 0$.

Situate sheet $a$ vertically above or below (depending on whether $a > 0$ or $a < 0$) the relevant edge of the auxiliary strip. Pick a place where sheet $a$ will not obscure any important existing intersections. Place another long, thin horizontal sheet so that its top edge is flush with the side of sheet $a$ opposite the auxiliary strip. Let $Y$ denote the intersection of the top edge of this last horizontal sheet and the zero crease.

Points $Z$ and $Y$ are on the zero crease and are separated by a distance of $xd + a$. We want to encode this distance as the separation between a new pair of long, thin horizontal strips in the same way that we were given $d$. To do this, take a pair of long, thin parallel rectangular strips so that their inner edges pass through $Z$ and $Y$. Use the diagonal reference to orient these strips at a 45° angle. Place a vertical strip across this pair of diagonal strips in a location where, looking horizontally from the intersections, you have a clean part of the square sheet in which to work. (If you have to go too far, it means that you did not start with a big enough square sheet.)

Finally, place a pair of thin, long horizontal strips so that their inner edges pass through the intersections of the inner edges of the pair of diagonal strips with the vertical strip. Take note of which horizontal strip traces back to point $Z$ and regard this strip as the one "locating zero".

## 7. THE CONSTRUCTION

As indicated in Section 5, to complete the construction, we iterate the procedure explained in previous section $n$ times. Begin the process by sandwiching the sheet representing the coefficient $a_n$ between two long, thin horizontal strips.

Now imagine sliding sheet $x$ left or right to vary $x$. As this is done, the entire contraption of horizontal, vertical, and diagonal sheets move in tandem according the prescription of alignments given for the construction. Undoubtedly, it would help to have a lot of people working as a team to do this! (A word of advice: put your fastest sprinters on the strips that encode $p(x)$.) The upshot is that the final pair of horizontal strips will widen and narrow representing the values of $p(x)$. When their inner edges touch, $x$ is a root of $p(x)$.

Our set up is slightly more amenable to computing positive real roots of polynomials, but this does not present any difficulty because one can always replace the polynomial $p(x)$ by $p(-x)$.

## 8. Returning to Origami Land

Given a polynomial, there are well-known bounds on the size of its real roots. This enables us to find sufficiently large, yet finite sheets to perform the massive procedure outlined in Sections 5-7.

Although our multifold uses several auxiliary parameters, the alignment $p(x) = 0$ depends only on the parameter $x$, so this is considered a one-parameter multifold. Also, by carefully following the procedure with all its recommendations, there will be no operational concerns such as having some important intersections or sheets blocked from access. Furthermore, the entire construction of $p(x)$ involves no unfolding, so we do not have to concern ourselves with the possibility of such unfolding causing a violation of condition 3 for the folding sequence of a multifold.

All that remains then, is to show how the whole contraption can be made from a single sheet of square origami paper!

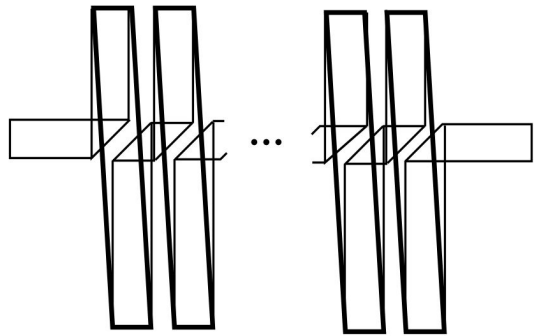

FIGURE 3. A schematic showing how to turn a strip into a rectangle. In actuality, all intersecting lines make 45° or 90° angles. In order to show all the folds clearly, we drew the folds askew to avoid overlaps. With minor alterations, the dimensions of the rectangle can be changed and the places where the strip enters and exits the rectangle can be relocated.

This can be done by connecting all sheets involved by an ultra thin strip of paper. When making these connections, attach the ultra thin strips to unimportant locations around the edges of the sheets used in the apparatus. Use plenty of slack in these connections so that the whole apparatus can be manipulated without obstructions. Now, purchase a really, really, really long, ultra thin strip of paper and fold it, using standard techniques (see Figure 3), into this massive contraption of connected sheets. Finally, if you are truly determined to use a square sheet of origami paper, purchase a gigantic square of origami paper with side length equal to the length of the super long ultra thin strip and fold it over and over along parallel creases until it has the exact dimensions of said super long ultra thin strip.

## 9. MULTIFOLDS AND BEYOND

The introduction of multifolds into origami constructibility fuels questions of practicality. Typically, multifolds are not easy to perform in the traditional way in which origami is practiced: by a lone folder. Indeed, the multifold introduced in this paper to construct $p(x)$ from $x$ would be impossible to fold by any normal human being! In view of this, should multifolds even be considered? We wish to make a few comments on this matter.

First, the concept of rolling a crease is not new or impractical. Folders roll through creases as a matter of routine. The second author recalls a lesson from Michael Lafosse that addressed how to fold a crease that passes precisely through a narrow tip. (For example, folding the lower left corner of a square to meet a horizontal crease over a crease that passes through the lower right corner...a maneuver often used to fold an equilateral triangle.) Michael demonstrated how one can achieve precision by pinching the crease at the tip in roughly the correct orientation and then rolling it to its final position before making it sharp.

Second, the one-parameter multifold is actually quite practical if one works in a team or builds a tool that helps by holding folds in place or keeps strategic parts of the paper taut. It would not be difficult, for instance, to build a tool that makes Lang's angle quintisection as simple as pulling on a loose end of paper. For an example of such a contraption that enables one to divide a sheet into $n$ equal sections, see Figure 4.

While on the question of what procedures are allowable in origami constructibility, we ask, can the world of origami constructible numbers be extended any further? Can non-algebraic numbers like pi be constructed as well?

Because alignments are specified by algebraic conditions, transcendental numbers would only be possible if the toolbox of allowable origami maneuvers is radically expanded.

To this end, one could introduce a mathematical model of an ideal sheet of paper that includes a description of how paper behaves in all three dimensions. Such a model would consider not only folds but also how paper curves in space. For example, one might consider continuous, path-length preserving maps from the square into space with the property that any region of the image that does not involve creases has zero curvature and satisfies certain physical conditions on tension. An allowable procedure would then be to specify a boundary condition in terms of existing references and extracting lengths between reference points in the resulting structure where possible. This would bring analysis into the world of origami constructible numbers and pi would likely be attainable.

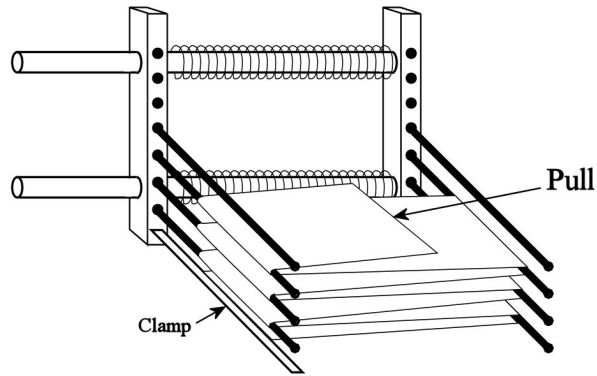

FIGURE 4. A custom tool that helps one make a one-parameter multifold. Pulling the free edge of a sheet woven around the cast iron metal rods against the tension of the two springs until the edge aligns with the topmost rod on the right allows one to divide the sheet into seven equal strips.